# MEASURES WITH ZEROS IN THE INVERSE OF THEIR MOMENT MATRIX

BY J. WILLIAM HELTON,[1] JEAN B. LASSERRE[2] AND MIHAI PUTINAR[3]

*University of California, San Diego, LAAS-CNRS and Institute of Mathematics, Toulouse and University of California, Santa Barbara*

We investigate and discuss when the inverse of a multivariate truncated moment matrix of a measure $\mu$ has zeros in some prescribed entries. We describe precisely which pattern of these zeroes corresponds to independence, namely, the measure having a product structure. A more refined finding is that the key factor forcing a zero entry in this inverse matrix is a certain *conditional triangularity property* of the orthogonal polynomials associated with $\mu$.

**1. Introduction.** It is well known that zeros in off-diagonal entries of the inverse $M^{-1}$ of a $n \times n$ covariance matrix $M$ identify pairs of random variables that have no *partial correlation* (and so are conditionally independent in case of normally distributed vectors); see, for example, Wittaker [7], Corollary 6.3.4. Allowing zeros in the off-diagonal entries of $M^{-1}$ is particularly useful for Bayesian estimation of regression models in statistics, and is called Bayesian *covariance selection*. Indeed, estimating a covariance matrix is a difficult problem for large number of variables, and exploiting sparsity in $M^{-1}$ may yield efficient methods for Graphical Gaussian Models (GGM). For more details, the interested reader is referred to Cripps, Carter and Kohn [3] and the many references therein.

The covariance matrix can be thought of as a matrix whose entries are second moments of a measure. This paper focuses on the truncated moment matrices, $M_d$, consisting of moments up to an order determined by $d$. First, we describe precisely the pattern of zeroes of $M_d^{-1}$ resulting from the measure having a product type structure. Next, we turn to the study

Received February 2007; accepted August 2008.
[1]Supported in part by the NSF and Ford Motor Co.
[2]Supported in part by the (French) ANR.
[3]Supported in part by the NSF.

of a particular entry of $M_d^{-1}$ being zero. We find that the key is what we call the *conditional triangularity* property of *orthonormal polynomials* (OP) up to degree $2d$, associated with the measure. To give the flavor of what this means, let, for instance, $\mu$ be the joint distribution $\mu$ of $n$ random variables $X = (X_1, \ldots, X_n)$, and let $\{p_\sigma\} \subset \mathbb{R}[X]$ be its associated family of orthonormal polynomials. When $(X_k)_{k \neq i,j}$ is fixed, they can be viewed as polynomials in $\mathbb{R}[X_i, X_j]$. If in doing so they exhibit a triangular structure [whence, the name conditional triangularity w.r.t. $(X_k)_{k \neq i,j}$], then entries of $M_d^{-1}$ at precise locations vanish. Conversely, if these precise entries of $M_d^{-1}$ vanish (robustly to perturbation), then the conditional triangularity w.r.t. $(X_k)_{k \neq i,j}$ holds. And so, for the covariance matrix case $d = 1$, this conditional triangularity property is equivalent to the zero partial correlation property well studied in statistics (whereas in general, we shall show that conditional independence is *not* detected by zeros in the inverse of the covariance matrix). Inverses of moment matrices naturally appear also in the recent work [2].

Interestingly, in a different direction, one may relate this issue with a constrained *matrix completion* problem. Namely, given that the entries of $M_d$ corresponding to marginals of the linear functional w.r.t. one variable at a time are fixed, complete the missing entries with values that make $M_d$ positive definite. This is a *constrained* matrix completion problem as one has to respect the moment matrix structure when filling up the missing entries. Usually, for the classical matrix completion problem with *no* constraint on $M$, the solution which maximizes an appropriate entropy gives *zeros* to entries of $M^{-1}$ corresponding to missing entries of $M$. But under the additional constraint of respecting the moment matrix structure, the maximum entropy solution does not always fill in $M_d^{-1}$ with zeros at the corresponding entries (as seen in examples by the authors). Therefore, any solution of this constrained matrix completion problem does not always maximize the entropy. Its "physical" or probabilistic interpretation is still to be understood.

We point out another accomplishment of this paper. More generally than working with a measure is working with a linear functional $\ell$ on the space of polynomials. One can consider moments with respect to $\ell$ and moment matrices. Our results hold at this level of generality.

**2. Notation and definitions.** For a real symmetric matrix $A \in \mathbb{R}^{n \times n}$, the notation $A \succ 0$ (resp. $A \succeq 0$) stands for $A$ positive definite (resp. semidefinite), and for a matrix $B$, let $B'$ or $B^T$ denote its transpose.

2.1. *Monomials, polynomials an moments.* We now discuss monomials at some length, since they are used in many ways, even to index the moment matrices which are the subject of this paper. Let $\mathbb{N}$ denote the nonnegative integers and $\mathbb{N}^n$ denote $n$ tuples of them and for $\alpha = (\alpha_1, \alpha_2 \cdots \alpha_n) \in \mathbb{N}^n$, define $|\alpha| := \sum_i \alpha_i$. The set $\mathbb{N}^n$ sits in one-to-one correspondence with the monomials via

$$\alpha \in \mathbb{N}^n \sim X^\alpha := X_1^{\alpha_1} X_2^{\alpha_2} \cdots X_n^{\alpha_n}.$$

Recall also the standard notation

$$\deg X^\alpha = |\alpha| = \alpha_1 + \cdots + \alpha_n.$$

By abuse of notation, we will freely interchange below $X^\alpha$ with $\alpha$, for instance, speaking about $\deg \alpha$ rather than $\deg X^\alpha$, and so on.

Let $\mathbb{R}[X]$ denote the ring of real polynomials in the variables $X_1, \ldots, X_n$ and let $\mathbb{R}_d[X] \subset \mathbb{R}[X]$ be the $\mathbb{R}$-vector space of polynomials of degree at most $d$. A polynomial $p \in \mathbb{R}[X]$ is a finite linear combination of monomials and it can be written

$$(2.1) \qquad p(X) = \sum_{\alpha \in \mathbb{N}^n} p_\alpha X^\alpha.$$

Let $y = (y_\alpha)_{\alpha \in \mathbb{N}^n}$ and define the linear functional $L_y : \mathbb{R}[X] \to \mathbb{R}$ first on monomials by $L_y(X^\alpha) = y_\alpha$ and then by linear extension to polynomials. That is,

$$(2.2) \qquad p \mapsto L_y(p) := \sum_{\alpha \in \mathbb{N}^n} p_\alpha y_\alpha,$$

whenever $(p_\alpha)$ are the coefficients of a polynomial as in (2.1). A linear functional $L_y$ on polynomials which is nonnegative (resp. positive) on all squares of polynomials [i.e., $L_y(p^2)$ is nonnegative] is what we call a *nonnegative (resp. positive) functional*.

The most prominent example is when the $y_\alpha$ are *moments of a measure* $\mu$ on $\mathbb{R}^n$, that is,

$$y_\alpha = \int_{\mathbb{R}^n} x^\alpha \mu(dx) \qquad \forall \alpha \in \mathbb{N}^n,$$

assuming of course that the signed measure $\mu$ decays fast enough at infinity, so that all monomials are integrable with respect to its total variation. Then

$$L_y(p) = \int_{\mathbb{R}^n} p(x) \, d\mu(x)$$

for every $p \in \mathbb{R}[X]$. Call $\mu$ a *representing measure* for $y$. For a positive measure $\mu$, the functional $L_y$ is nonnegative, however, in the converse direction, a linear functional $L_y$ on polynomials being nonnegative on squares of polynomials is not equivalent to there existing a positive measure $\mu$. This is closely related to Hilbert's 17th Problem and its progeny, focusing on positive polynomials not always being a sum of squares.

As the reader will see much of what we do here holds for positive linear functionals, no measure is required. To state our results, we must introduce finite moment matrices. Their entries are indexed by monomials and so we must describe orders on monomials.

2.2. *Orders on monomials.* Define the *fully graded partial order* (*FG order*) "≤", on monomials, or equivalently, at the level of multi-indices, define $\gamma \leq \alpha$ for $\gamma, \alpha \in \mathbb{N}^n$ iff $\gamma_j \leq \alpha_j$ for all $j = 1, \ldots, g$. Important to us is

$$\alpha \leq \beta \quad \text{iff } X^\alpha \text{ divides } X^\beta.$$

Define the *graded lexicographic order* (*GLex*) "$<_{\mathrm{gl}}$", on monomials, or equivalently, for $\gamma, \alpha \in \mathbb{N}^n$ first by using $\deg \gamma \leq \deg \alpha$ to create a partial order. Next refine this to a total order by breaking ties in two monomials $X^{m_1}$, $X^{m_2}$ of the same degree $|m_1| = |m_2|$, as would a dictionary with $X_1 = a, X_2 = b, \ldots$. For example, the monomials in two variables $X_1, X_2$ of degree $\leq 2$ listed in GLex order are

$$1, X_1, X_2, X_1^2, X_1 X_2, X_2^2.$$

Beware $\gamma <_{\mathrm{gl}} \alpha$ does not imply $\gamma \leq \alpha$; for example, $(1,1,3) <_{\mathrm{gl}} (1,3,1)$, but $\leq$ fails. However, $\beta \leq \alpha$ and $\beta \neq \alpha$ imply $\beta <_{\mathrm{gl}} \alpha$.

It is convenient to list all monomials as an infinite vector $v_\infty(X) := (X^\alpha)_{\alpha \in \mathbb{N}^n}$, where the entries are listed in GLex order, henceforth called the *tautological vector*; $v_d(X) = (X^\alpha)_{|\alpha| \leq d} \in \mathbb{R}^{s(d)}$ denotes the finite vector consisting of the part of $v_\infty(X)$ containing exactly the degree $\leq d$ monomials. Using this notation, we can write polynomials as

(2.3) $$p(X) = \langle \mathbf{p}, v_d(X) \rangle$$

for some real vector $\mathbf{p} = (p_\alpha)$, where the latter is the standard nondegenerate pairing between $\mathbb{R}^{s(d)}$ and $\mathbb{R}^{s(d)} \otimes_\mathbb{R} \mathbb{R}[X]$.

2.3. *Moment matrix.* Given a sequence $y = (y_\alpha)_{\alpha \in \mathbb{N}^n}$, the *moment matrix* $M_d(y)$ associated with $y$ has its rows and columns indexed by $\alpha$, $|\alpha| \leq d$, where the $\alpha$ are listed in GLex order, and

$$M_d(y)(\alpha, \beta) := L_y(X^\alpha X^\beta) = y_{\alpha+\beta} \quad \forall \alpha, \beta \in \mathbb{N}^n \text{ with } |\alpha|, |\beta| \leq d.$$

For example, $M_2(y)$ is

$$M_2(y): \begin{array}{c} \\ 1 \to \\ X_1 \to \\ X_2 \to \\ X_1^2 \to \\ X_1 X_2 \to \\ X_2^2 \to \end{array} \begin{pmatrix} 1 & X_1 & X_2 & X_1^2 & X_1 X_2 & X_2^2 \\ 1 & y_{10} & y_{01} & y_{20} & y_{11} & y_{02} \\ y_{10} & y_{20} & y_{11} & y_{30} & y_{21} & y_{12} \\ y_{01} & y_{11} & y_{02} & y_{21} & y_{12} & y_{03} \\ y_{20} & y_{30} & y_{21} & y_{40} & y_{31} & y_{22} \\ y_{11} & y_{21} & y_{12} & y_{31} & y_{22} & y_{13} \\ y_{02} & y_{12} & y_{03} & y_{22} & y_{13} & y_{04} \end{pmatrix}.$$

Note that the functional $L_y$ produces for every $d \geq 0$ a positive semidefinite moment matrix $M_d(y)$ if and only if

$$L_y(p^2) \geq 0 \qquad \forall p \in \mathbb{R}[X],$$

that is, iff $L_y$ is a nonnegative functional. The associated matrices $M_d(y)$ are positive definite, for all $d \geq 0$ if and only if $L_y(p^2) = 0$ implies $p = 0$.

**3. Measures of product form.** Now we can describe one pursuit of this paper. Given a sequence $y = (y_\alpha)$ indexed by $\alpha, |\alpha| \leq 2d$, we investigate some properties of the inverse $M_d(y)^{-1}$ of a positive definite moment matrix $M_d(y)$ when entries of the latter satisfy a product form property.

DEFINITION 1. We say that the moment matrix $M_d(y) \succ 0$ has the *product form property*, if

$$(3.1) \qquad L_y(X^\alpha) = \prod_{i=1}^n L_y(X_i^{\alpha_i}) \qquad \forall \alpha \in \mathbb{N}^n, |\alpha| \leq 2d,$$

or equivalently, we say the positive linear functional $L_y$ has *the independence property*. If $M_d(y)^{-1}(\alpha, \beta) = 0$ for every $y$ such that $M_d(y) \succ 0$ has the product form property, then we say the pair $(\alpha, \beta)$ is a *congenital zero for d-moments*.

For example, if $y$ consists of moments of a product measure $\mu = \prod_i^n \mu_i(x_i)$, with $\int \mu_i = 1$, then (3.1) corresponds to the fact

$$L_y(X^\alpha) = \int_{\mathbb{R}^n} X^\alpha \mu = \prod_i^n \int_{\mathbb{R}} X_i^{\alpha_i} \mu_i.$$

For random variables, congenital zeroes can be thought of as those zeroes in $M_d(y)^{-1}$ due to independence.

We now can give the flavor of our main results.

THEOREM 3.1. *The pair $(\alpha, \beta) \in \mathbb{N}^n \times \mathbb{N}^n$ is a congenital zero for the d-moment problem if and only if the least common multiple $X^\eta$ of $X^\alpha$ and $X^\beta$ has degree bigger than $d$.*

The above result can be conveniently rephrased in terms of the max operation defined for $\alpha, \beta \in \mathbb{N}^n$ by

$$\max(\alpha, \beta) := (\max(\alpha_j, \beta_j))_{j=1,\ldots,n}.$$

Set $\eta := \max(\alpha, \beta)$. Simple observations about this operation are as follows:

1. $X^\eta$ is the least common multiple of $X^\alpha$ and $X^\beta$,

2. $X^\alpha$ divides $X^\beta$ iff $X^\beta = X^\eta$,
   $X^\beta$ divides $X^\alpha$ iff $X^\alpha = X^\eta$,
3. $|\eta| = \sum_{j=1}^{n} \max(\alpha_j, \beta_j)$.

Thus, Theorem 3.1 asserts that the entry $(\alpha, \beta)$ does not correspond to a congenital zero in the matrix $M_d(y)^{-1}$ if and only if $|\max(\alpha, \beta)| \leq d$.

Later in Theorem 5.1 we show that this LCM (least common multiple) characterization of zeros in $M_d^{-1}(y)$ is equivalent to a highly triangular structure of orthonormal polynomials associated with the positive functional $L_y$.

EXAMPLE. In the case of $M_2^{-1}(y)$ in two variables $X_1, X_2$ we indicate below which entries $M_2(y)^{-1}(\alpha, \beta)$, with $|\beta| = 2$, are congenital zeroes. These $(\alpha, \beta)$ index the last three columns of $M_2(y)^{-1}$ and are

$$
\begin{array}{rcccc}
 & & X_1^2 & X_1 X_2 & X_2^2 \\
1 & \to & * & * & * \\
X_1 & \to & * & * & 0 \\
X_2 & \to & 0 & * & * \\
X_1^2 & \to & * & 0 & 0 \\
X_1 X_2 & \to & 0 & * & 0 \\
X_2^2 & \to & 0 & 0 & *.
\end{array}
$$

Here $*$ means that the corresponding entry can be different from zero. Note each $*$ corresponds to $X^\alpha$ failing to divide $X^\beta$.

The proof relies on properties of orthogonal polynomials, so we begin by explaining in some detail the framework.

3.1. *Orthonormal polynomials.* A functional analytic viewpoint to polynomials is expeditious, so we begin with that. Let $s(d) := \binom{n+d}{d}$ be the dimension of vector space $\mathbb{R}_d[X]$. Let $\langle \cdot, \cdot \rangle_{\mathbb{R}^{s(d)}}$ denote the standard inner product on $\mathbb{R}^{s(d)}$. Let $f, h \in \mathbb{R}[X]$ be the polynomials $f(X) = \sum_{|\alpha|=0}^{s(d)} f_\alpha X^\alpha$ and $h(X) = \sum_{|\alpha|=0}^{s(d)} h_\alpha X^\alpha$. Then,

$$\langle f, h \rangle_y := \langle \mathbf{f}, M_d(y)\mathbf{h} \rangle_{\mathbb{R}^{s(d)}} = L_y(f(X)h(X))$$

defines a *scalar product* in $\mathbb{R}_d[X]$, provided $M_d(y)$ is positive definite.

With a given $y = (y_\alpha)$ such that $M_d(y) \succ 0$, one may associate a unique family $(p_\alpha)_{|\alpha|=0}^{s(d)}$ of *orthonormal polynomials*. That is, the $p_\alpha$'s satisfy

(3.2) $\begin{cases} p_\alpha \in \text{lin.span}\{X^\beta; \beta \leq_{\text{gl}} \alpha\}, \\ \langle p_\alpha, p_\beta \rangle_y = \delta_{\alpha\beta}, & |\alpha|, |\beta| \leq d, \\ \langle p_\alpha, X^\beta \rangle_y = 0, & \text{if } \beta <_{\text{gl}} \alpha, \langle p_\alpha, X^\alpha \rangle_y > 0. \end{cases}$

Note $\langle p_\alpha, X^\beta \rangle_y = 0$, if $\beta \leq \alpha$ and $\alpha \neq \beta$, since the latter implies $\beta <_{\mathrm{gl}} \alpha$.

Existence and uniqueness of such a family is guaranteed by the Gram–Schmidt orthonormalization process following the GLex order on the monomials, and by the positivity of the moment (covariance) matrix; see, for instance, [1], Theorem 3.1.11, page 68.

*Computation.* Although not needed for the rest of the present article, a determinantal formula for the orthogonal polynomials is within reach, with a proof very similar to the classical one in the one variable case. The reader can omit this subsection without loss of continuity.

Suppose that we want to compute the orthonormal polynomials $p_\sigma$ for some index $\sigma$. Then proceed as follows: build up the sub-moment matrix $M^{(\sigma)}(y)$ with columns indexed by all monomials $\beta \leq_{\mathrm{gl}} \sigma$, and rows indexed by all monomials $\alpha <_{\mathrm{gl}} \sigma$. Hence, $M^{(\sigma)}(y)$ has one row less than columns. Next, complete $M^{(\sigma)}(y)$ with an additional last row described by $[M^{(\sigma)}(y)]_{\sigma,\beta} = X^\beta$, for all $\beta \leq_{\mathrm{gl}} \sigma$. Then up to a normalizing constant, $p_\sigma$ is nothing less than $\det(M^\sigma(y))$.

To see this, let $\gamma <_{\mathrm{gl}} \sigma$. Then

$$\langle X^\gamma, p_\sigma \rangle_y = L_y(X^\gamma \det(M^{(\sigma)}(y))) = \det(B^\sigma)(y),$$

where the matrix $B^\sigma(y)$ is the same as $M^\sigma(y)$ except for the last row which is now the vector $(L_y(X^{\gamma+\alpha}))_{\alpha \leq_{\mathrm{gl}} \sigma}$ already present in one of the rows above. Therefore, $\det(B^\sigma)(y) = 0$. For instance, with $n = 2$ and the ordering $X_1 < X_2$, let $\sigma := (0, 1)$. Then $\langle X_1, p_\sigma \rangle_y = 0$ because

$$L_y\left(X_1 \det \begin{bmatrix} 1 & y_{10} & y_{01} \\ y_{10} & y_{20} & y_{11} \\ 1 & X_1 & X_2 \end{bmatrix}\right) = L_y\left(\det \begin{bmatrix} 1 & y_{10} & y_{01} \\ y_{10} & y_{20} & y_{11} \\ X_1 & X_1^2 & X_1 X_2 \end{bmatrix}\right)$$

$$= \det \begin{bmatrix} 1 & y_{10} & y_{01} \\ y_{10} & y_{20} & y_{11} \\ y_{10} & y_{20} & y_{11} \end{bmatrix} = 0,$$

and similarly, $\langle 1, p_\sigma \rangle_y = 0$.

Next, writing its coefficient $\rho_{\sigma\beta}$ is just (again up to a normalizing constant) the cofactor of the element $[M^\sigma(y)]_{1,\beta}$ in the square matrix $M^\sigma(y)$ with rows and columns both indexed with $\alpha \leq_{\mathrm{gl}} \sigma$.

*Further properties.* Now we give further properties of the orthonormal polynomials. Consider first one variable polynomials. The orthogonal polynomials $p_k$ have, by their very definition, a "triangular" form, namely

$$p_k(X_1) := \sum_{\ell \leq k} \rho_{k\ell} X_1^\ell.$$

The orthonormal polynomials inherit the product form property of $M_d(y)$, assuming that the latter holds. Namely, each orthonormal polynomial $p_\alpha$ is a product

$$p_\alpha(X) = p_{\alpha_1}(X_1) p_{\alpha_2}(X_2) \cdots p_{\alpha_n}(X_n) \tag{3.3}$$

of orthogonal polynomials $p_{\alpha_j}(X_j)$ in one dimension. Indeed, by the product property,

$$\langle p_{\alpha_1}(X_1) p_{\alpha_2}(X_2) \cdots p_{\alpha_n}(X_n), X^\beta \rangle_y = \prod_{j=1}^n \langle p_{\alpha_j}(X_j), X_j^{\beta_j} \rangle_y,$$

whence the product of single variable orthogonal polynomials satisfies all requirements listed in (3.2).

"Triangularity" in one variable and the product form property (3.3) forces $p_\alpha$ to have what we call a *fully triangular* form:

$$p_\alpha(X) := \sum_{\gamma \leq \alpha} \rho_{\alpha\gamma} X^\gamma, \qquad |\alpha| \leq d. \tag{3.4}$$

Also note that *for any $\gamma \leq \alpha$ there exists a positive functional $L_y$ of product type making $\rho_{\alpha\gamma}$ not zero.*

To exhibit such a functional, we will use a particular property of coefficients of Laguerre polynomials. Given $\sigma \in \mathbb{N}^n$, consider the product measure $\mu(dx) = \prod_{i=1}^n \mu_i(dx_i)$ on the positive orthant $\mathbb{R}_+^n$, with $\mu_i(dx_i) = e^{-x_i} x_i^{\sigma_i} \, dx_i$.

The univariate Laguerre polynomials

$$x_i \mapsto L_k^{(\sigma_i)}(x_i) = \frac{e^{x_i} x_i^{-\sigma_i}}{k!} \frac{d^k}{dx_i^k}(e^{-x_i} x_i^{n+\sigma_i}) = \sum_{j=0}^k \binom{k+\sigma}{k-j} \frac{(-x_i)^j}{j!}$$

are orthogonal with respect to the measure $\mu_i$ on the semi-axis $\mathbb{R}_+$; see, for instance, [4]. Observe that the degree of the coefficient of $x_i^j$ with respect to the variable $\sigma_i$ is precisely $k - j$.

The orthogonal polynomials associated with the product measure $\mu$ and its associated positive functional

$$L_y(p) = \int_{\mathbb{R}_+^n} p(x_1, \ldots, x_n) e^{-x_1 - \cdots - x_n} x_1^{\sigma_1} \cdots x_n^{\sigma_n} \, dx_1 \cdots dx_n, \qquad p \in \mathbb{R}[X],$$

are the (Laguerre$_\sigma$) polynomials $\Lambda_\alpha^\sigma(X) = \prod_{j=1}^n L_{\alpha_j}^{(\sigma_j)}(X_j)$.

We formalize a simple observation as a lemma because we use it later.

LEMMA 3.2. *The coefficients $\rho_{\alpha,\beta}(\sigma)$ in the decomposition*

$$\Lambda_\alpha^\sigma(X) = \sum_{\beta \leq \alpha} \rho_{\alpha,\beta}(\sigma) X^\beta$$

*of Laguerre$_\sigma$ polynomials are themselves polynomials in $\sigma = (\sigma_1, \ldots, \sigma_n)$, viewed as independent variables, and the multi-degree of $\rho_{\alpha,\beta}(\sigma)$ is $\alpha - \beta$.*

Note that for an appropriate choice of the parameters $\sigma_j, 1 \leq j \leq n$, the coefficients $\rho_{\alpha,\beta}(\sigma)$ in the decomposition

$$\Lambda_\alpha^\sigma(X) = \sum_{\beta \leq \alpha} \rho_{\alpha,\beta}(\sigma) X^\beta$$

are linearly independent over the rational field, and hence nonnull. To prove this, evaluate $\sigma$ on an $n$-tuple of algebraically independent transcendental real numbers over the rational field.

Note that we have used very little of the specific properties of the Laguerre polynomials. The mere precise *polynomial* dependence of their coefficients with respect to the variable $\sigma$ was sufficient for the proof of the above lemma.

3.2. *Proof of Theorem* 3.1. Let $L_y$ be a linear functional for which $M_d(y) \succ 0$ and let $(p_\alpha)$ denote the family of orthogonal polynomials with respect to $L_y$. Orthogonality in (3.2) for expansions (3.4) reads

$$\delta_{\alpha\beta} = \langle p_\alpha, p_\beta \rangle_y = \sum_{\gamma \leq \alpha, \sigma \leq \beta} \rho_{\alpha\gamma} \rho_{\beta\sigma} \langle X^\gamma, X^\sigma \rangle_y.$$

In matrix notation this is just

$$I = DM_d(y)D^T,$$

where $D$ is the matrix $D = (\rho_{\alpha\gamma})_{|\alpha|,|\gamma| \leq d}$. Its columns are indexed (as before) by monomials arranged in GLex order, and likewise for its rows. That $\rho_{\alpha\gamma} = 0$ if $\gamma \not\leq \alpha$ implies that $\rho_{\alpha\gamma} = 0$ if $\alpha <_{\text{gl}} \gamma$, which says precisely that $D$ is lower triangular. Moreover, its diagonal entries $\rho_{\beta\beta}$ are not 0, since $p_\beta$ must have $X^\beta$ as its highest order term. Because of this and triangularity, $D$ is invertible. Write

$$M_d(y) = D^{-1}(D^T)^{-1} \quad \text{and} \quad M_d(y)^{-1} = D^T D.$$

Our goal is to determine which entries of $M_d(y)^{-1}$ are forced to be zero and we proceed by writing the formula $Z := M_d(y)^{-1} = D^T D$ as

$$
\begin{aligned}
z_{\alpha\beta} &= \sum_{|\gamma| \leq d} \rho_{\gamma\alpha} \rho_{\gamma\beta} = \sum_{\beta \leq \gamma, \alpha \leq \gamma, |\gamma| \leq d} \rho_{\gamma\alpha} \rho_{\gamma\beta} \\
&= \sum_{\max(\alpha,\beta) \leq \gamma, |\gamma| \leq d} \rho_{\gamma\alpha} \rho_{\gamma\beta}.
\end{aligned}
$$
(3.5)

We emphasize (since it arises later) that this uses only the full triangularity of the orthogonal polynomials rather than that they are products of one variable polynomials. If full triangularity were replaced by triangularity w.r. to $<_{\text{gl}}$, then the first two equalities in (3.5) would be the same except that $\beta \leq \gamma, \alpha \leq \gamma, |\gamma| \leq d$ would be replaced by $\beta \leq_{\text{gl}} \gamma, \alpha \leq_{\text{gl}} \gamma, |\gamma| \leq d$.

To continue with our proof, consider $(\alpha, \beta)$ and set $\eta := \max(\alpha, \beta)$. If $|\max(\alpha, \beta)| > d$, then $z_{\alpha,\beta} = 0$, since the sum in equation (3.5) is empty. This is the forward side of Theorem 3.1.

Conversely, consider the product measure $\mu$ on the positive orthant $\mathbb{R}_+^n$ whose associated orthogonal polynomials are the Laguerre$_\sigma$ polynomials $\{\Lambda_\alpha^\sigma\}$ of Lemma 3.2. When $|\max(\alpha, \beta)| \leq d$ the entry $z_{\alpha\beta}$ is a sum of one or more products $\rho_{\gamma\alpha}\rho_{\gamma\beta}$ and so is a polynomial in $\sigma$. If this polynomial is not identically zero, then some value of $\sigma$ makes $z_{\alpha\beta} \neq 0$, so $(\alpha, \beta)$ is not a congenital zero. Now we set out to show that $z_{\alpha\beta}$ (as a polynomial in $\sigma$) is not identically 0.

Lemma 3.2 tells us each product $\rho_{\gamma,\alpha}\rho_{\gamma,\beta}$ is a polynomial whose multi-degree in $\sigma$ is exactly $2\gamma - \alpha - \beta$. The multi-index $\gamma$ is subject to the constraints $\max(\alpha, \beta) \leq \gamma$ and $|\gamma| \leq d$. We fix an index, say, $j = 1$, and choose

$$\hat{\gamma} = \max(\alpha, \beta) + (d - |\max(\alpha, \beta)|, 0, \ldots, 0).$$

Note the product $\rho_{\hat{\gamma}\alpha}\rho_{\hat{\gamma}\beta}$ is included in the sum (3.5) for $z_{\alpha,\beta}$ and it is a polynomial of degree $2d - 2|\max(\alpha, \beta)| + 2\max(\alpha_1, \beta_1) - \alpha_1 - \beta_1$ with respect to $\sigma_1$. By the extremality of our construction of $\hat{\gamma}$, every other term $\rho_{\gamma\alpha}\rho_{\gamma\beta}$ in $z_{\alpha\beta}$ will have smaller $\sigma_1$ degree. Hence, $\rho_{\hat{\gamma}\alpha}\rho_{\hat{\gamma}\beta}$ cannot be canceled, proving that $z_{\alpha\beta}$ is not the zero polynomial.

**4. Partial independence.** In this section we consider the case where only a partial independence property holds. We decompose the variables into disjoint sets $X = (X(1), \ldots, X(k))$, where $X(1) = (X(1)_1, \ldots, X(1)_{d_1})$, and so on. Note that the lexicographic order on $\mathbb{N}^n$ respects the grouping of variables, in the sense

$$(\alpha_1, \ldots, \alpha_k) <_{\text{gl}} (\beta_1, \ldots, \beta_k)$$

if and only if, either $\alpha_1 <_{\text{gl}} \beta_1$, or, if $\alpha_1 = \beta_1$, then, either $\alpha_2 <_{\text{gl}} \beta_2$, and so on.

The linear functional $L_y$ is said to satisfy a *partial independence property* (w.r.t. the fixed grouping of variables), if

$$L_y(p_1(X(1)) \cdots p_k(X(k))) = \prod_{j=1}^{k} L_y(p_j(X(j)),$$

where $p_j$ is a polynomial in the variables from the set $X(j)$, respectively.

In this context we still, analogously to Definition 1, use the term congenital zeros in connection with inverses of moment matrices $M(y) \succ 0$ corresponding to $L_y$ having the partial independence property. Now we state the natural generalization of Theorem 3.1 to partial independence.

THEOREM 4.1. *Let $L_y$ be a positive functional satisfying a partial independence property with respect to the groups of variables $X = (X(1), \ldots, X(k))$. Let $\alpha = (\alpha_1, \ldots, \alpha_k), \beta = (\beta_1, \ldots, \beta_k)$ be two multi-indices decomposed according to the fixed groups of variables, and satisfying $|\alpha|, |\beta| \leq d$.*

*Then the $(\alpha, \beta)$-entry in the matrix $M_d(y)^{-1}$ is congenitally zero if and only if, for every $\gamma = (\gamma_1, \ldots, \gamma_k)$ satisfying $\gamma_j \geq_{\mathrm{gl}} \alpha_j, \beta_j, 1 \leq j \leq k$, we have $|\gamma| > d$.*

The structure behind this is just the analog of what we used before. Denote by $\deg_{X(j)} Q(X)$ the degree of a polynomial $Q$ in the variables $X(j)$. Assuming that $L_y$ is a positive functional, one can associate in a unique way the orthogonal polynomials $p_\alpha, \alpha \in \mathbb{N}^n$. Let $\alpha = (\alpha_1, \ldots, \alpha_k)$ be a multi-index decomposed with respect to the groupings $X(1), \ldots, X(k)$. Then, the uniqueness property of the orthonormal polynomials implies

$$p_\alpha(X) = p_{\alpha_1}(X(1)) \cdots p_{\alpha_k}(X(k)),$$

where $p_{\alpha_j}(X(j))$ are orthonormal polynomials depending solely on $X(j)$, and arranged in lexicographic order within this group of variables.

With this background, the proof of Theorem 4.1 repeats that of Theorem 3.1, with only the observation that

$$p_\alpha(X) = \sum_{\substack{1 \leq j \leq k \\ \gamma_j \leq_{\mathrm{gl}} \alpha_j}} c_{\alpha,\gamma} X^\gamma.$$

**5. Full triangularity.** A second look at the proof of Theorem 3.1 reveals that the only property of the multivariate orthogonal polynomials we have used was the full triangularity form (3.4). In this section we provide an example of a nonproduct measure which has orthogonal polynomials in full triangular form, and, on the other hand, we prove that the zero pattern appearing in our main result, in the inverse moment matrices $M_r(y)^{-1}, r \leq d$, implies the full triangular form of the associated orthogonal polynomials. Therefore, *zeros in the inverse $M_d^{-1}$ are coming from a certain triangularity property of orthogonal polynomials rather than from a product form of $M_d$*.

EXAMPLE 1. We work in two real variables $(x, y)$, with the measure $d\mu = (1 - x^2 - y^2)^t \, dx \, dy$, restricted to the unit disk $x^2 + y^2 < 1$, where $t > -1$ is a parameter.

Let $P_k(u; s)$ denote the orthonormalized Jacobi polynomials, that is, the univariate orthogonal polynomials on the segment $[-1, 1]$, with respect to the measure $(1 - u^2)^s \, du$, with $s > -1$.

According to [6], Example 1, Chapter X, the orthonormal polynomials associated to the measure $d\mu$ on the unit disk are

$$Q_{m+n,n}(x, y) = P_m(x, t + n + 1/2)(1 - x^2)^{n/2} P_n\left(\frac{y}{\sqrt{1 - x^2}}; t\right)$$

and

$$Q_{n,m+n}(x,y) = P_m(y, t+n+1/2)(1-y^2)^{n/2} P_n\left(\frac{x}{\sqrt{1-y^2}}; t\right).$$

Observe that these polynomials have full triangular form, yet the generating measure is not a product of measures.

THEOREM 5.1 (Full triangularity theorem). *Let $y = (y_\alpha)_{\alpha \in \mathbb{N}^n}$ be a multi-sequence, such that the associated moment matrices $M_d(y)$ are positive definite, where $d$ is a fixed positive integer. Then the following holds.*

*For every $r \leq d$, the $(\alpha, \beta)$-entry in $M_r(y)^{-1}$ is $0$ whenever $|\max(\alpha, \beta)| > r$ if and only if the associated orthogonal polynomials $P_\alpha, |\alpha| \leq d$, have full triangular form.*

PROOF. The proof of the forward side is exactly the same as in the proof of Theorem 3.1. The converse proof begins by expanding each orthogonal polynomial $p_\alpha$ as

(5.1) $$p_\alpha(X) := \sum_{\gamma \leq_{\text{gl}} \alpha} \rho_{\alpha\gamma} X^\gamma, \qquad |\alpha| \leq d,$$

where we emphasize that $\gamma \leq_{\text{gl}} \alpha$ is used as opposed to (3.4) and we want to prove that $\beta \leq_{\text{gl}} \alpha$ and $\beta \geq \alpha$ imply $\rho_{\alpha\beta} = 0$. Note that when $|\alpha| = d$ the inequality

(5.2) $$\beta \geq \alpha \text{ is equivalent to } |\max(\alpha, \beta)| > d.$$

The first step is to use the zero locations of $M_d(y)^{-1}$ to prove that

(5.3) $$\rho_{\alpha\beta} = 0 \qquad \text{if } |\max(\alpha, \beta)| > d,$$

only for $|\alpha| = d$. Once this is established, then we apply the same argument to prove (5.2) for $|\alpha| = d'$ where $d'$ is smaller.

If $n = 1$, there is nothing to prove. Assume $n > 1$ and decompose $\mathbb{N}^n$ as $\mathbb{N} \times \mathbb{N}^{n-1}$. The corresponding indices will be denoted by $(i, \alpha), i \in \mathbb{N}, \alpha \in \mathbb{N}^{n-1}$. We shall prove (5.3) by descending induction on the graded lexicographic order applied to the index $(i, \alpha)$ the following statement:

*Let $i + |\alpha| = d$ and assume that $(j, \beta) <_{\text{gl}} (i, \alpha)$ and $(j, \beta) \geq (i, \alpha)$. Then $\rho_{(i,\alpha),(j,\beta)} = 0$.*

The precise statement we shall use is equivalent [because of (5.2)] to the following.

*The induction hypothesis*: *Suppose that*

(5.4) $$\rho_{(i',\alpha'),(j,\beta)} = 0 \qquad \text{if } \max(i', j) + |\max(\alpha', \beta)| > d,$$

*holds for all indices $(i', \alpha') >_{\text{gl}} (i, \alpha)$, with $i' + |\alpha'| = i + |\alpha| = d$.*

We want to prove $\rho_{(i,\alpha),(j,\beta)} = 0$. Since we shall proceed by induction from the top, we assume that $i + |\alpha| = d$ that $(j,\beta) \geq (i,\alpha)$ that $(j,\beta) <_{\text{gl}} (i,\alpha)$ and let $(i,\alpha)$ denote the largest such w.r.t. $<_{\text{gl}}$ order. Clearly, $i = d$. There is only one corresponding term in the graded lexicographic sequence of indices of length less than or equal to $d$, namely, $(d,0)$. We shall prove that, for every $\beta \in \mathbb{N}^{n-1}, (j,\beta) \leq_{\text{gl}} (d,0), |\beta| > 0$, we have

$$\rho_{(d,0),(j,\beta)} = 0.$$

Since the corresponding entry in $M_d(y)^{-1}$, denoted henceforth as before by $z_{**}$,

$$z_{(d,0),(j,\beta)} = 0,$$

is zero by assumption, and because

$$z_{(d,0),(j,\beta)} = \rho_{(d,0),(d,0)} \rho_{(d,0),(j,\beta)},$$

we obtain $\rho_{(d,0),(j,\beta)} = 0$.

Now we turn to proving the main induction step. Assuming (5.4), we want to prove $\rho_{(i,\alpha),(j,\beta)} = 0$. Let $(j,\beta) \leq_{\text{gl}} (i,\alpha)$ subject to the condition $\max(i,j) + |\max(\alpha,\beta)| > d$. Then by hypothesis $0 = z_{(i,\alpha),(j,\beta)}$, so the GLex version of expansion (3.5) gives

$$0 = \rho_{(i,\alpha),(i,\alpha)} \rho_{(i,\alpha),(j,\beta)} + \sum_{i'=i, \alpha' >_{\text{gl}} \alpha, i'+|\alpha'|=d} \rho_{(i',\alpha'),(i,\alpha)} \rho_{(i',\alpha'),(j,\beta)}$$
$$+ \sum_{i'>i, i'+|\alpha'|=d} \rho_{(i',\alpha'),(i,\alpha)} \rho_{(i',\alpha'),(j,\beta)}.$$

We will prove that the two summations above are zero.

Indeed, if $i' > i$, then $\max(i,i') + |\max(\alpha',\alpha)| > i + |\alpha| = d$, and the induction hypothesis implies $\rho_{(i',\alpha'),(i,\alpha)} = 0$, which eliminates the second sum.

Assume $i' = i$, so that $|\alpha| = |\alpha'|$. Then if $\max(\alpha,\alpha')$ equals to either $\alpha$ or $\alpha'$, we get $\alpha = \alpha'$, but this cannot be, since $\alpha' >_{\text{gl}} \alpha$. Thus, $i + |\max(\alpha,\alpha')| > d$ and the induction hypothesis yields in this case $\rho_{(i',\alpha'),(i,\alpha)} = 0$. Which eliminates the first sum.

In conclusion,

$$\rho_{(i,\alpha),(i,\alpha)} \rho_{(i,\alpha),(j,\beta)} = 0,$$

which implies

$$\rho_{(i,\alpha),(j,\beta)} = 0,$$

as desired.

Our induction hypothesis is valid and we initialized it successfully, so we obtain the working hypothesis:

$$\rho_{(i,\alpha),(j,\beta)} = 0,$$

whenever $(j, \beta) <_{\text{gl}} (i, \alpha)$, $|(i, \alpha)| = d$, and $\max((j, \beta), (i, \alpha)) > d$. By (5.2), this is full triangularity under the assumption $i + |\alpha| = d$. □

In contrast to the above full triangularity criteria, the product decomposition of a potential truncated moment sequence $(y_\alpha)_{\alpha_i \leq d}$ can be decided by elementary linear algebra. Assume for simplicity that $n = 2$, and write the corresponding indices as $\alpha = (i, j) \in \mathbb{N}^2$. Then there are numerical sequences $(u_i)_{i \leq d}$ and $(v_j)_{j \leq d}$ with the property

$$y_{(i,j)} = u_i v_j, \qquad 0 \leq i, j \leq d,$$

if and only if

$$\operatorname{rank}(y_{(i,j)})_{i,j=0}^{d} \leq 1.$$

A similar rank condition, for a corresponding multilinear map, can be deduced for arbitrary $n$.

**6. Conditional triangularity.** The aim of this section is to extend the full triangularity theorem of Section 5, to a more general context.

6.1. *Conditional triangularity.* In this new setting we consider two tuples of variables

$$X = (x_1, \ldots, x_n), \qquad Y = (y_1, \ldots, y_m),$$

and we will impose on the concatenated tuple $(X, Y)$ the full triangularity only with respect to the set of variables $Y$. The conclusion is as expected: this assumption will reflect the appearance of some zeros in the inverse $M_d^{-1}$ of the associated truncated moment matrix. The proof below is quite similar to that of Theorem 5.1 and we indicate only sufficiently many details to make clear the differences.

We denote the set of indices by $(\alpha, \beta)$, with $\alpha \in \mathbb{N}^n, \beta \in \mathbb{N}^m$, equipped with the graded lexicographic order "$<_{\text{gl}}$". In addition, the set of indices $\beta \in \mathbb{N}^m$, which refers to the set of variables $Y$, is also equipped with the full graded order "$\leq$".

Let $y = (y_\alpha)_{\alpha \in \mathbb{N}^n}$ be a multi-sequence, such that the associated moment matrices $M_d(y)$ are positive definite, where $d$ is a positive integer. As before, we denote

$$p_{(\alpha,\beta)}(X, Y) = \sum_{(\alpha,\beta) \geq_{\text{gl}} (\alpha',\beta')} \rho_{(\alpha,\beta),(\alpha',\beta')} X^{\alpha'} Y^{\beta'}$$

the associated orthogonal polynomials, and by

$$z_{(\alpha,\beta),(\alpha',\beta')} = \sum_{(\gamma,\sigma) \geq_{\text{gl}} (\alpha,\beta),(\alpha',\beta')} \rho_{(\gamma,\sigma),(\alpha,\beta)} \rho_{(\gamma,\sigma),(\alpha',\beta')},$$

the entries in $M_d(y)^{-1}$.

DEFINITION 2 (*Conditional triangularity*). The orthonormal polynomials $p_{\alpha,\beta} \in \mathbb{R}[X,Y]$, with $|\alpha + \beta| \leq 2d$, satisfy the *conditional triangularity with respect to* $X$, if when $X$ is fixed and considered as a parameter, the resulting family denoted $\{p_{\alpha,\beta}|X\} \subset \mathbb{R}[Y]$ is in full triangular form with respect to the $Y$ variables. More precisely, the following $(O)_d$ condition below holds:

$$(O)_d : [(\alpha', \beta') \leq_{\mathrm{gl}} (\alpha, \beta), |(\alpha, \beta)| \leq d, \text{ and } \beta' \geq \beta] \Rightarrow \rho_{(\alpha,\beta),(\alpha',\beta')} = 0.$$

Next, for a fixed degree $d \geq 1$, we will have to consider the following *zero in the inverse condition*.

DEFINITION 3 [*Zero in the inverse condition* $(V)_d$]. Assume the degree $d$ is fixed. Let $(\alpha', \beta') \leq_{\mathrm{gl}} (\alpha, \beta)$ with $|(\alpha, \beta)| \leq d$ be arbitrary.
If $|(\gamma, \sigma)| > d$ whenever $(\gamma, \sigma) \geq_{\mathrm{gl}} (\alpha, \beta), (\alpha', \beta')$ and $\sigma \geq \max(\beta, \beta')$, then $z_{(\alpha,\beta),(\alpha',\beta')} = 0$.

The main result of this paper asserts that both conditions $(V)_d$ and $(O)_d$ are in fact equivalent.

THEOREM 6.1 (Conditional triangularity). *Let $y = (y_\alpha)_{\alpha \in \mathbb{N}^n}$ be a multi-sequence and let $d$ be an integer such that the associated moment matrices $M_d(y)$ are positive definite.*
*Then the zero in the inverse condition $(V)_r$, $r \leq d$, holds if and only if $(O)_d$ holds, that is, if and only if the orthonormal polynomials satisfy the conditional triangularity with respect to $X$.*

PROOF. Clearly, from its definition, $(O)_d$ implies $(O)_r$ for all $r \leq d$. One direction is obvious: Let $r \leq d$ be fixed, arbitrary. If condition $(O)_r$ holds, then a pair $((\alpha', \beta') \leq_{\mathrm{gl}} (\alpha, \beta))$ subject to the assumptions in $(V)_r$ will leave not a single term in the sum giving $z_{(\alpha,\beta),(\alpha',\beta')}$.

Conversely, assume that $(V)_d$ holds. We will prove the vanishing statement $(O)_d$ by descending induction with respect to the graded lexicographical order. To this aim, we label all indices $(\alpha, \beta), |(\alpha, \beta)| = d$ in decreasing graded lexicographic order:

$$(\alpha_0, \beta_0) >_{\mathrm{gl}} (\alpha_1, \beta_1) >_{\mathrm{gl}} (\alpha_2, \beta_2) \ldots.$$

In particular, $d = |\alpha_0| \geq |\alpha_1| \geq \cdots$ and $0 = |\beta_0| \leq |\beta_1| \leq \cdots$.

To initialize the induction, consider $(\alpha', \beta') \leq_{\mathrm{gl}} (\alpha_0, \beta_0) = (\alpha_0, 0)$, with $\beta' \geq 0$, that is, $\beta' \neq 0$. Then

$$0 = z_{(\alpha_0,\beta_0),(\alpha',\beta')} = \rho_{(\alpha_0,\beta_0),(\alpha_0,\beta_0)} \rho_{(\alpha_0,\beta_0),(\alpha',\beta')}.$$

Since the leading coefficient $\rho_{(\alpha_0,\beta_0),(\alpha_0,\beta_0)}$ in the orthogonal polynomial is nonzero, we infer

$$[(\alpha',\beta') <_{\rm gl} (\alpha_0,\beta_0), \beta' \geq \beta_0] \implies \rho_{(\alpha_0,\beta_0),(\alpha',\beta')} = 0,$$

which is exactly condition $(O)_d$ applied to this particular choice of indices.

Assume that $(O)_d$ holds for all $(\alpha_j,\beta_j), j < k$. Let $(\alpha',\beta') <_{\rm gl} (\alpha_k,\beta_k)$ with $\beta' \geq \beta_k$, that is, $|\max(\beta',\beta_k)| > |\beta_k|$. In view of $(V_d)$,

$$z_{(\alpha_k,\beta_k),(\alpha',\beta')} = \rho_{(\alpha_k,\beta_k),(\alpha_k,\beta_k)}\rho_{(\alpha_k,\beta_k),(\alpha',\beta')} + \sum_{j=0}^{k-1} \rho_{(\alpha_j,\beta_j),(\alpha_k,\beta_k)}\rho_{(\alpha_j,\beta_j),(\alpha',\beta')}.$$

Note that the induction hypothesis implies that, for every $0 \leq j \leq k-1$, at least one factor $\rho_{(\alpha_j,\beta_j),(\alpha_k,\beta_k)}$ or $\rho_{(\alpha_j,\beta_j),(\alpha',\beta')}$ vanishes. Thus,

$$0 = z_{(\alpha_k,\beta_k),(\alpha',\beta')} = \rho_{(\alpha_k,\beta_k),(\alpha_k,\beta_k)}\rho_{(\alpha_k,\beta_k),(\alpha',\beta')},$$

whence $\rho_{(\alpha_k,\beta_k),(\alpha',\beta')} = 0$.

Once we have exhausted by the above induction all indices of length $d$, we proceed similarly to those of length $d-1$, using now the zero in the inverse property $(V)_{d-1}$, and so on. $\square$

Remark that the ordering $(X,Y)$ with the tuple of full triangular variables $Y$ on the second entry is important. A low degree example will be considered in the last section.

We call Theorem 6.1 the conditional triangularity theorem because when the variables $X$ are *fixed* (and so can be considered as parameters), then the orthogonal polynomials now considered as elements of $\mathbb{R}[Y]$ are in full triangular form, rephrased as *in triangular form conditional to $X$ is fixed*.

6.2. *The link with partial correlation.* Let us specialize to the case $d = 1$. Assume that the underlying joint distribution on the random vector $X = (X_1,\ldots,X_n)$ is *centered*, that is, $\int X_i\,d\mu = 0$ for all $i = 1,\ldots,n$; then $M_d$ reads

$$M_d = \left[\begin{array}{c|c} 1 & 0 \\ \hline 0 & R \end{array}\right] \quad \text{and} \quad M_d^{-1} = \left[\begin{array}{c|c} 1 & 0 \\ \hline 0 & R^{-1} \end{array}\right],$$

where $R$ is just the usual *covariance* matrix. Partitioning the random vector as $(Y, X_i, X_j)$ with $Y = (X_k)_{k \neq i,j}$ yields

$$R = \begin{bmatrix} \operatorname{var}(Y) & \operatorname{cov}(Y,X_i) & \operatorname{cov}(Y,X_j) \\ \operatorname{cov}(Y,X_i) & \operatorname{var}(X_i) & \operatorname{cov}(X_i,X_j) \\ \operatorname{cov}(Y,X_j) & \operatorname{cov}(X_i,X_j) & \operatorname{var}(X_j) \end{bmatrix},$$

where var and cov have obvious meanings. The *partial covariance* of $X_i$ and $X_j$ given $Y$, denoted $\text{cov}(X_i, X_j | Y)$ in Wittaker ([7], page 135) satisfies:

(6.1)  $\text{cov}(X_i, X_j | Y) := \text{cov}(X_i, X_j) - \text{cov}(Y, X_i) \text{var}(Y)^{-1} \text{cov}(Y, X_j).$

After scaling $R^{-1}$ to have a unit diagonal, the *partial correlation* between $X_i$ and $X_j$ (partialled on $Y$) is the negative of $\text{cov}(X_i, X_j | Y)$, and as already mentioned, $R^{-1}(i,j) = 0$ if and only if $X_i$ and $X_j$ have *zero partial correlation*, that is, $\text{cov}(X_i, X_j | Y) = 0$. See, for example, Wittaker [7], Corollaries 5.8.2 and 5.8.4.

COROLLARY 6.2. *Let $d = 1$. Then $R^{-1}(i,j) = 0$ if and only if the orthonormal polynomials of degree up to 2, associated with $M_1$, satisfy the conditional triangularity with respect to $X = (X_k)_{k \neq i,j}$.*

PROOF.  To recast the problem in the framework of Section 6.1, let $Y = (X_i, X_j)$ and rename $X := (X_k)_{k \neq i,j}$. In view of Definition 3 with $d = 1$, we only need consider pairs $(\alpha', \beta') \leq_{\text{gl}} (\alpha, \beta)$ with $\alpha = \alpha' = 0$ and $\beta' = (0,1)$, $\beta = (1,0)$. But then $\sigma \geq \max[\beta, \beta'] = (1,1)$ implies $|(\gamma, \sigma)| \geq 2 > d$, and so as $R^{-1}(i,j) = 0$, the zero in the inverse condition $(V)_d$ holds. Equivalently, by Theorem 6.1, $(O)_d$ holds. □

Corollary 6.2 states that the pair $(X_i, X_j)$ has zero partial correlation if and only if the orthonormal polynomials up to degree 2 satisfy the conditional triangularity with respect to $X = (X_k)_{k \neq i,j}$. That is, partial correlation and conditional triangularity are *equivalent*.

EXAMPLE 2.  Let $d = 1$, and consider the case of three random variables $(X, Y, Z)$ with (centered) joint distribution $\mu$. Then suppose that the orthonormal polynomials up to degree $d = 1$ satisfy the conditional triangularity property $(V)_1$ w.r.t. $X$. That is, $p_{000} = 1$ and

$$p_{100} = \alpha_1 + \beta_1 X,$$
$$p_{010} = \alpha_2 + \beta_2 X + \gamma_2 Y,$$
$$p_{001} = \alpha_3 + \beta_3 X + \gamma_3 Z,$$

for some coefficients $(\alpha_i \beta_i \gamma_i)$. Notice that because of $(O)_1$, we cannot have a linear term in $Y$ in $p_{001}$. Orthogonality yields that $\langle X^\gamma, p_\alpha \rangle = 0$ for all $\gamma <_{\text{gl}} \alpha$, that is, with **E** being the expectation w.r.t. $\mu$,

$$\beta_2 \mathbf{E}(X^2) + \gamma_2 \mathbf{E}(XY) = 0,$$
$$\beta_3 \mathbf{E}(X^2) + \gamma_3 \mathbf{E}(XZ) = 0,$$
$$\beta_3 \mathbf{E}(XY) + \gamma_3 \mathbf{E}(YZ) = 0.$$

Stating that the determinant of the last two linear equations in $(\beta_3, \gamma_3)$ is zero yields

$$\mathbf{E}(YZ) - \mathbf{E}(X,Y)\mathbf{E}(X^2)^{-1}\mathbf{E}(X,Z) = 0,$$

which is just (6.1), that is, the zero partial correlation condition, up to a multiplicative constant.

This immediately raises two questions:

(i) What are the distributions for which the orthonormal polynomials up to degree 2 satisfy the conditional triangularity with respect to a given pair $(X_i, X_j)$?

(ii) Among such distributions, what are those for which conditional independence with respect to $X = (X_{k \neq i,j})$ also holds?

An answer to the latter would characterize distributions for which zero partial correlation imply conditional independence (like for the normal distribution).

6.3. *Conditional independence and zeros in the inverse.* We have already mentioned that, in general, conditional independence is *not* detected from zero entries in the inverse of $R^{-1}$ (equivalently, $M_1^{-1}$), except for the normal joint distribution, a common assumption in Graphical Gaussian Models. Therefore, a natural question of potential interest is to search for conditions on when conditional independence in the non-Gaussian case is related to the *zero in the inverse* property $(V)_d$, or equivalently, the conditional triangularity $(O)_d$, not only for $d = 1$, but also for $d > 1$.

A rather negative result in this direction is as follows. Let $d$ be fixed, arbitrary, and let $M_d = (y_{ijk}) \succ 0$ be the moment matrix of an arbitrary joint distribution $\mu$ of three random variables $(X, Y_1, Y_2)$ on $\mathbb{R}$. As we are considering only finitely many moments (up to order $2d$), by Tchakaloff's theorem, there exists a measure $\varphi$ finitely supported on, say $s$, points $(x^{(l)}, y_1^{(l)}, y_2^{(l)}) \subset \mathbb{R}^3$ (with associated probabilities $\{p_l\}$), $l = 1, \ldots, s$, and whose moments up to order $2d$ match those of $\mu$; see, for example, Reznick [5], Theorem 7.18.

Let us define a sequence $\{\varphi_t\}$ of probability measures as follows. Perturb each point $(x^{(l)}, y_1^{(l)}, y_2^{(l)})$ to $(x^{(l)} + \epsilon(t,l), y_1^{(l)}, y_2^{(l)})$, $l = 1, \ldots, s$, in such a way that no two points $x^{(l)} + \epsilon(t,l)$ are the same, and keep the same weights $\{p_l\}$. It is clear that $\varphi_t$ satisfies

$$1 = \text{Prob}[Y = (y_1^{(l)}, y_2^{(l)}) | X = x^{(l)} + \epsilon(t,l)]$$
$$= \text{Prob}[Y_1 = y_1^{(l)} | X = x^{(l)} + \epsilon(t,l)]\text{Prob}[Y_2 = y_2^{(l)} | X = x^{(l)} + \epsilon(t,l)]$$

for all $l = 1, \ldots, s$. That is, conditional to $X$, the variables $Y_1$ and $Y_2$ are independent. Take a sequence with $\epsilon(t, l) \to 0$ for all $l$, as $t \to \infty$, and consider the moment matrix $M_d^{(t)}$ associated with $\varphi_t$. Clearly, as $t \to \infty$,

$$\int X^i Y_1^j Y_2^k \, d\varphi_t \to \int X^i Y_1^j Y_2^k \, d\mu = y_{ijk} \qquad \forall i, j, k : i + j + k \leq 2d,$$

that is, $M_d^{(t)} \to M_d$.

Therefore, if the zero in the inverse property $(V)_d$ does not hold for $M_d$, then by a simple continuity argument, it cannot hold for any $M_d^{(t)}$ with sufficiently large $t$, and still the conditional independence property holds for each $\varphi_t$. One has just shown, that for every fixed $d$, one may easily construct examples of measures with the conditional independence property, which violate the zero in the inverse property $(V)_d$.

**Acknowledgments.** This work was initiated at Banff Research Station (Canada) during a workshop on Positive Polynomials, and completed at the Institute for Mathematics and its Applications (IMA, Minneapolis). All three authors wish to thank both institutes for their excellent working conditions and financial support.

J. W. HELTON  
DEPARTMENT OF MATHEMATICS  
UNIVERSITY OF CALIFORNIA AT SAN DIEGO  
LA JOLLA, CALIFORNIA 92093  
USA  
E-MAIL: helton@math.ucsd.edu  

J. B. LASSERRE  
LAAS-CNRS AND INSTITUTE OF MATHEMATICS  
7 AVENUE DU COLONEL ROCHE  
31077 TOULOUSE CÉDEX 4  
FRANCE  
E-MAIL: lasserre@laas.fr

M. Putinar
Department of Mathematics
University of California at Santa Barbara
Santa Barbara, California 93106
USA
E-mail: mputinar@math.ucsb.edu